\documentclass[preprint]{elsarticle} 

\usepackage{amscd,amsxtra,amsbsy,amsfonts,amssymb,amstext,amsxtra,amsmath,amsthm}
\usepackage{a4,array,delarray,textcomp}
\usepackage{graphicx}
\usepackage{float}
\usepackage{pspicture,pst-all,multido}
\usepackage[latin1]{inputenc}
\usepackage{ngerman,url}
\usepackage{euler}

\newcommand{\C}{\mathbb{C}}
\newcommand{\R}{\mathbb{R}}
\newcommand{\N}{\mathbb{N}}

\restylefloat{figure}

\newcounter{envcount}%

  {\vspace{\bigskipamount}}

  {\vspace{\bigskipamount}}

  {\vspace{\bigskipamount}}

  {\vspace{\bigskipamount}}

  {\vspace{\bigskipamount}}

  {\vspace{\bigskipamount}}

  {\vspace{\bigskipamount}}

  {\vspace{\bigskipamount}}

  {\vspace{\bigskipamount}}

  {\vspace{\bigskipamount}}

  {\vspace{\bigskipamount}}

  {\vspace{\bigskipamount}}

  {\vspace{\bigskipamount}}    

  {\vspace{\bigskipamount}}    

  {\vspace{\bigskipamount}}

  {\vspace{\bigskipamount}}

  {\vspace{\bigskipamount}}

  {\vspace{\bigskipamount}}

  {\vspace{\bigskipamount}}

{\vspace{\bigskipamount}\refstepcounter{envcount}\textbf{(\theenvcount)\enspace Lokale Multiplizit\"at und Kardinalit\"at der Einpunkturbilder.}}%
  {\vspace{\bigskipamount}}

{\vspace{\bigskipamount}\refstepcounter{envcount}\textbf{(\theenvcount)\enspace Separable $L^2$-R\"aume.}}%
  {\vspace{\bigskipamount}}

{\vspace{\bigskipamount}\refstepcounter{envcount}\textbf{(\theenvcount)\enspace Spektraldarstellung mit Vielfachheiten eines normaler Operators.}}%
  {\vspace{\bigskipamount}}

{\vspace{\bigskipamount}\refstepcounter{envcount}\textbf{(\theenvcount)\enspace Multiplikationsoperator.}}%
  {\vspace{\bigskipamount}}

{\vspace{\bigskipamount}\refstepcounter{envcount}\textbf{(\theenvcount)\enspace Funktionenalgebra Isomorphismus.}\itshape}%
  {\vspace{\bigskipamount}\upshape}

{\vspace{\bigskipamount}\refstepcounter{envcount}\textbf{(\theenvcount)\enspace  Von Neumann Modell.} \itshape}%
  {\vspace{\bigskipamount}\upshape}

{\vspace{\bigskipamount}\refstepcounter{envcount}\textbf{(\theenvcount)\enspace   Zur Eindeutigkeit des Von Neumann Modells.} \itshape}%
  {\vspace{\bigskipamount}\upshape}

{\vspace{\bigskipamount}\refstepcounter{envcount}\textbf{(\theenvcount)\enspace Spektraldarstellung mit Vielfachheiten aus dem Von Neumann Modell.} \itshape}%
  {\vspace{\bigskipamount}\upshape}

{\vspace{\bigskipamount}\refstepcounter{envcount}\textbf{(\theenvcount)\enspace Ma\ss algebra.}}%
  {\vspace{\bigskipamount}}

{\vspace{\bigskipamount}\refstepcounter{envcount}\textbf{(\theenvcount)\enspace Analytischer Messraum.}}%
  {\vspace{\bigskipamount}}

{\vspace{\bigskipamount}\refstepcounter{envcount}\textbf{(\theenvcount)\enspace Standardma\ss raum.}}%
  {\vspace{\bigskipamount}}

{\vspace{\bigskipamount}\refstepcounter{envcount}\textbf{(\theenvcount)\enspace Orthogonale Zerlegung.}}%
  {\vspace{\bigskipamount}}

  {\vspace{\bigskipamount}}

{\vspace{\bigskipamount}\refstepcounter{envcount}\textbf{(\theenvcount)\enspace Vorbereitung.}}%
  {\vspace{\bigskipamount}}

{\vspace{\bigskipamount}\refstepcounter{envcount}\textbf{(\theenvcount)\enspace Zusammenfassung.}}%
  {\vspace{\bigskipamount}}

{\vspace{\bigskipamount}\refstepcounter{envcount}\textbf{(\theenvcount)\enspace Definition.}}%
  {\vspace{\bigskipamount}}

{\vspace{\bigskipamount}\enspace \textbf{Definition.}}%
  {\vspace{\bigskipamount}}

{\vspace{\bigskipamount}\enspace \fbox{Üb}}%
  {\vspace{\bigskipamount}}

{\vspace{\bigskipamount}\refstepcounter{envcount}\textbf{(\theenvcount)\enspace Satz.}\itshape}%
  {\vspace{\bigskipamount}\upshape}

\newenvironment{The}%
{\vspace{\bigskipamount}\refstepcounter{envcount}\textbf{(\theenvcount)\enspace Theorem.}\itshape}%
  {\vspace{\bigskipamount}\upshape}

\newenvironment{Pro}%
{\vspace{\bigskipamount}\refstepcounter{envcount}\textbf{(\theenvcount)\enspace Proposition.}\itshape}%
  {\vspace{\bigskipamount}\upshape}

{\vspace{\bigskipamount}\refstepcounter{envcount}\textbf{(\theenvcount)\enspace Erinnerung.}\itshape}%
  {\vspace{\bigskipamount}\upshape}

{\vspace{\bigskipamount}\refstepcounter{envcount}\textbf{(\theenvcount)\enspace Korollar.}\itshape}%
  {\vspace{\bigskipamount}\upshape}

{\vspace{\bigskipamount}\refstepcounter{envcount}\textbf{(\theenvcount)\enspace Corollary.}\itshape}%
  {\vspace{\bigskipamount}\upshape}

 \newenvironment{Exa}%
{\vspace{\bigskipamount}\refstepcounter{envcount}\textbf{(\theenvcount)\enspace Example.}\itshape}%
  {\vspace{\bigskipamount}\upshape}

{\vspace{\bigskipamount}\refstepcounter{envcount}\textbf{(\theenvcount)\enspace Definition und Satz.}}%
  {\vspace{\bigskipamount}}

{\vspace{\bigskipamount}\refstepcounter{envcount}\textbf{(\theenvcount)\enspace Satz und Definition.}}%
  {\vspace{\bigskipamount}}

{\vspace{\bigskipamount}\refstepcounter{envcount}\textbf{(\theenvcount)\enspace Definition und Korollar.}}%
  {\vspace{\bigskipamount}}

{\vspace{\bigskipamount}\refstepcounter{envcount}\textbf{(\theenvcount)\enspace Definition und Lemma.}}%
  {\vspace{\bigskipamount}}

\newenvironment{Lem}%
{\vspace{\bigskipamount}\refstepcounter{envcount}\textbf{(\theenvcount)\enspace Lemma.}\itshape}%
  {\vspace{\bigskipamount}\upshape}

{\vspace{\bigskipamount}\refstepcounter{envcount}\textbf{(\theenvcount)\enspace Lemma und Definition.}\itshape}%
  {\vspace{\bigskipamount}\upshape}

{\vspace{\bigskipamount}\refstepcounter{envcount}\textbf{(\theenvcount)\enspace Folgerung.}}%
  {\vspace{\bigskipamount}}

{\vspace{\bigskipamount}\refstepcounter{envcount}\textbf{(\theenvcount)\enspace Folgerungen.}}%
  {\vspace{\bigskipamount}}

{\vspace{\bigskipamount}\refstepcounter{envcount}\textbf{(\theenvcount)\enspace Zusatz.}}%
  {\vspace{\bigskipamount}}

{\vspace{\bigskipamount}\refstepcounter{envcount}\textbf{(\theenvcount)\enspace Zusammenfassung.}}%
  {\vspace{\bigskipamount}}

{\vspace{\bigskipamount}\refstepcounter{envcount}\textbf{(\theenvcount)\enspace Beispiel.}}%
  {\vspace{\bigskipamount}}

{\vspace{\bigskipamount}\refstepcounter{envcount}\textbf{(\theenvcount)\enspace Beispiele.}}%
  {\vspace{\bigskipamount}}

{\vspace{\bigskipamount}\refstepcounter{envcount}\textbf{(\theenvcount)\enspace Beispiel}}%
  {\vspace{\bigskipamount}}

{\vspace{\bigskipamount}\refstepcounter{envcount}\textbf{(\theenvcount)\enspace Bemerkung.}}%
  {\vspace{\bigskipamount}}

{\vspace{\bigskipamount}\refstepcounter{envcount}\textbf{(\theenvcount)\enspace Bemerkung}}%
  {\vspace{\bigskipamount}}

{\vspace{\bigskipamount}\refstepcounter{envcount}\textbf{(\theenvcount)\enspace Bemerkungen.}}%
  {\vspace{\bigskipamount}}

{\vspace{\bigskipamount}\refstepcounter{envcount}\textbf{(\theenvcount)\enspace}}%
  {\vspace{\bigskipamount}}

{\vspace{\bigskipamount}\refstepcounter{envcount}\textbf{(\theenvcount)\enspace Schreibweisen und Bezeichnungen}}%
  {\vspace{\bigskipamount}}

\newenvironment{PO}%
{\vspace{\bigskipamount}\refstepcounter{envcount}\textbf{(\theenvcount)\enspace Proof of}}%
  {\vspace{\bigskipamount}}

\theoremstyle{definition}
\swapnumbers

\theoremstyle{remark}

\begin{document}
\title{Cyclicity for Unbounded Multiplication Operators in $L^p$- and $C_0$\,-Spaces}

\author{Sebastian Zaigler} % \fnref{fn 1}
\ead{sebastian@zaigler.com}

\author{Domenico P.L. Castrigiano}
\ead{castrig@ma.tum.de}

\address{Technische Universit\"at M\"unchen, Fakult\"at f\"ur Mathematik, M\"unchen, Germany} 
%\fntext[fn 1]{Corresponding author}

\begin{abstract}

For every, possibly unbounded,  multiplication operator in $L^p$-space, $p\in\,]0,\infty[$, on finite separable measure space we show that multicyclicity, multi-$*$-cyclicity, and multiplicity coincide. This result includes and generalizes Bram's much cited theorem  from 1955 on bounded $*$-cyclic normal operators. It also includes as a core result cyclicity of the multiplication operator  $M_z$ by the complex variable $z$  in $L^p(\mu)$ for every Borel measure $\mu$ on $\C$. The concise proof is based in part  on the result that the function $e^{-\left|z\right|^2}$ is a $*$-cyclic vector for $M_z$ in $C_0(\C)$ and further  in $L^p(\mu)$. We characterize topologically those locally compact sets $X\subset \C$, for which $M_z$ in $C_0(X)$ is cyclic.

\end{abstract}


\maketitle

\section*{Introduction}

In 1955 Bram \cite{Br} proves his well-known  and much cited theorem that a bounded $*$-cyclic normal operator is cyclic. It is also well-known that, as a consequence, a normal operator is cyclic if and only if it has multiplicity one or, equivalently, if it is simple. 
In 2009 Nagy \cite{Na1}   tackles the generalization of Bram's result to unbounded normal operators. \\
\hspace*{6mm}  Due to the spectral theorem  the question actually  concerns  multiplication operators in $L^2(\mu)$ for  finite Borel measures on $\C$ with possibly unbounded support. We extend the frame to general (unbounded) multiplication operators in $L^p$-spaces for $p\in\,]0,\infty[$ on finite separable measure spaces.
 We prove that multicyclicity, multi-$*$-cyclicity, and multiplicity coincide for those operators.
\\
\hspace*{6mm} This result  includes   cyclicity of the multiplication operator $M_z$ by $z$ in $L^p(\mu)$ for any finite Borel measure $\mu$ on $\C$, which in turn  is a main step in the proof of the above result. It is not even obvious how to prove $*$-cyclicity of $M_z$ in the unbounded case,
and there are several futile attempts   in the literatur
concerning  the Hilbert space case. \\
\hspace*{6mm}  Let us rapidly recall the case of bounded $M_z$ in $L^2(\mu)$.  Here supp$(\mu)$ is compact. Then the set $\Pi(z,\overline{z})$ of polynomials in $z$ and $\bar{z}$ is dense in  $L^2(\mu)$, since $\Pi(z,\overline{z})$ is dense in $C(\operatorname{supp}(\mu))$ by the theorem of Stone/Weierstra\ss . Therefore, if $\Pi(z,\bar{z})$ is contained in the closure of the polynomials $\Pi(z)$, the latter are dense in $L^2(\mu)$,  i.e., $\overline{\Pi(z)}=L^2(\mu)$, which means that the constant $1_\C$ is a cyclic vector for $M_z$. Actually, still due to the boundedness of $M_z$, it suffices to show that the function $\bar{z}$ is element of the closure of $\Pi(z)$. Bram \cite{Br} solves this approximation problem decomposing $\C$ into the union of an increasing sequence of $\alpha$-sets and a $\mu$-null set.  An $\alpha$-set is a compact subset of $\C$  such that every continuous function on it can be approximated uniformly by polynomials in $z$.  By Lavrentev's theorem  (see, e.g. \cite{Car,Ha})   
the $\alpha$-sets are just the compact
subsets of $\C$ with empty interior and connected complement. \\
\hspace*{6mm} In the unbounded case this way has to be modified, mostly by two reasons.
First, due to unboundedness the support of $\mu$ is not  compact and  the polynomials  in $z$ and $\bar{z}$ are not bounded on $\operatorname{supp}(\mu)$.
Secondly,  $\Pi(z,\overline{z})$  need not be dense in $L^2(\mu)$ (see e.g. Hamburger's example  in Simon \cite[example 1.3]{Si}).  In \cite{Na1} Nagy generalizes Bram's decomposition of $\C$ for any (non-compact) $\operatorname{supp}(\mu)$.   This is an important result. We have considerably simplified its proof. By this \cite{Na1}  succeeds in showing that $\bar{z}$ is in the closure of $\Pi(z)$ in $L^p(\sigma)$ for some finite Borel measure $\sigma$ equivalent to $\mu$. The ensuing tacit assumption by Nagy that $z\,\overline{\Pi(z)}\subset L^p(\sigma)$ however, as we will explain below, definitely restricts the proof to the case of bounded $M_z$, thus missing the aim.  \\
\hspace*{6mm}
 As already implied,  for the unbounded case, in the literature  there seems even to exist no explicit proof   for $*$-cyclicity of $M_z$ in $L^2(\mu)$.  However, Agricola/Friedrich  \cite[sec. 3]{AF}  show that the functions $p\operatorname{e}^{-|x|^2}$,  $p$ polynomial on $\R^d$,  are dense in $C_0(\R^d)$ with respect to uniform convergence. 
 In particular this means that the function $e^{-|z|^2}$ is a $*$-cyclic vector  for $M_z$ in $C_0(\C)$. As a ready consequence,  $e^{-|z|^2}$ is $*$-cyclic for $M_z$ in $L^p(\mu)$. Then we proceed similarly to \cite{Na1},  but showing at once by an induction argument that the whole of $\Pi(z,\bar{z})$ lies in the closure of $\Pi(z)$. We like to remark that we present in (\ref{ISCVFCO})  a short classical proof  of the density result  of Agricola/Friedrich  (which is central in \cite{AF}) and that  we apply successfully the same method   for the proof of other results on cyclicity.
  \\
 \hspace*{6mm}
 Another important ingredient   is the  Rohlin  decomposition  of a measurable function which we apply to unbounded functions on  finite separable measure spaces.\\
  \hspace*{6mm}
We get started on the multiplication operator $M_z$ in $C_0(\C)$ and extend to $M_z$ in   $C_0(X)$ for  locally compact $X\subset\C$. We find that $M_z$ is $*$-cyclic  and  describe  topologically those $X$, for which $M_z$ is cyclic.\\
 \hspace*{6mm}
Finally, it is worth   mentioning  that the results on ($*$)-cyclicity  for the most part are obtained by  polynomial approximation, thus contributing to this field.  We shall give some examples.

\section*{Main results}

Let $p\in\,]0,\infty[$ and let $(\Omega,\mathcal{A},\mu)$ be a measure space. For a measurable function $\varphi: \Omega \to \C$ let $M_\varphi$ denote the multiplication operator in $L^p(\mu)$ given by $M_\varphi f:= \varphi f$ with domain $\mathcal{D}(M_\varphi):=\{f\in L^p(\mu):\,\varphi f \in L^p(\mu)\}$. We will deal with  separable $L^p$-spaces. Therefore it is no restriction to assume that $\mu$ is finite and that the measure algebra is separable. It is well-known that $M_\varphi$ is closed and,  if $p=2$, normal. Moreover, $M_\varphi$ is bounded iff $\varphi$ is $\mu$-essentially bounded.\\
\hspace*{6mm}
A set $Z\subset L^p(\mu)$ is called \textbf{cyclic} for $M_\varphi$ 
if   $p(\varphi)f\in L^p(\mu)$ for all $p\in\Pi(z)$, $f\in Z$ and if $$\Pi(M_\varphi)Z:=\{p(\varphi)f: p\in\Pi(z), f\in Z\}$$ is dense in $L^p(\mu)$. If there is no finite cyclic set  the \textbf{multicyclicity} mc$(M_\varphi)$ is  set $\infty$.  Otherwise  it is defined as the smallest number of elements of a cyclic set. $M_\varphi$ is called cyclic if mc$(M_\varphi)=1$. Similarly, a set $Z\subset L^p(\mu)$ is called \textbf{$*$-cyclic} for $M_\varphi$ 
if   $p(\varphi,\overline{\varphi})f\in L^p(\mu)$ for all $p\in\Pi(z,\bar{z})$, $f\in Z$ and if $$\Pi(M_\varphi, M_{\overline{\varphi}})Z:=\{p(\varphi,\overline{\varphi})f: p\in\Pi(z,\bar{z}), f\in Z\}$$ is dense in $L^p(\mu)$. The \textbf{multi-$*$-cyclicity} mc*$(M_\varphi)$ is defined analogously and $M_\varphi$ is called $*$-cyclic if mc*$(M_\varphi)=1$.\\
\hspace*{6mm}  We choose to define  multiplicity  by means of the Rohlin decomposition $(\pi,\nu)$ of $(\varphi,\mu)$. Let us briefly explain this decomposition. See also Seid  \cite[remark 3.4]{Se}.  There is a measure algebra isomorphism from $(\Omega,\mathcal{A},\mu)$ onto $([0,1]\times \C, \mathcal{B}, \nu)$ with $\mathcal{B}$ the Borel sets and $\nu$ a finite measure. The latter satisfies 
$$\nu=\lambda\otimes \mu_c +\sum_{n\in\N}\delta_{1/n}\otimes \mu_n,$$ 
where $\lambda$ denotes the Lebesgue measure on $[0,1]$, $\delta_{1/n}$ is the point measure at $1/n$, and $\mu_c$, $\mu_n$ are Borel measures on $\C$ with $\mu_{n+1} \ll \mu_n$ for $n\in\N$.  Moreover, by this measure algebra isomorphism, $M_\varphi$ in $L^p(\mu)$ is isomorphic with $M_\pi$ in $L^p(\nu)$ with $\pi(t,z):=z$.
Since the measures $\lambda\otimes \mu_c$, $\delta_1\otimes \mu_1$, $\delta_{1/2}\otimes\mu_2,\dots$  are mutually orthogonal, $L^p(\nu)$ and $M_\pi$  are identified with the $p$-direct sums

\begin{equation*}
%L^p(\nu)\simeq 
L^p(\lambda\otimes \mu_c) \oplus \bigoplus_{n\in\N} L^p(\mu_n), \quad M_\pi\oplus \bigoplus_{n\in \N} M_z.
\end{equation*}

Then $M_z$ in $L^p(\mu_n)$, $n\in\N$, is cyclic, whereas  $M_\pi$ on a subspace $\{1_Sf:f\in L^p(\lambda\otimes\mu_c)\}$ with $S\in \mathcal{B}$ is cyclic only if  $\lambda\otimes\mu_c(S)=0$. Hence, in view of $\mu_{n+1} \ll \mu_n$ for $n\in\N$,   the \textbf{multiplicity} of $M_\varphi$ is defined as mp$(M_\varphi) :=\sup\{n\in\N:\mu_n\ne 0\}$ if $\mu_c=0$ and $\infty$ otherwise.\\
\hspace*{6mm}  Note that for  $\mu_c\ne 0$  the normal operator  $M_\pi$ in $L^2(\lambda\otimes \mu_c)$ is  Hilbert space isomorphic  with the countably infinite orthogonal sum of copies of $M_z$ in $L^2(\mu_c)$. It follows that the above definition of multiplicity  coincides in the case $p=2$ with the usual multiplicity mp$(T)$ for a normal operator. Obviously this is not in contrast to the fact that due to the spectral theorem a normal operator $T$ in separable Hilbert space 
 is isomorphic with $M_\pi$  in $L^2(\nu)$ for a Rohlin decomposition with $\mu_c=0$. \\
\hspace*{6mm}  In case that $(\Omega,\mathcal{A})$ is a standard measurable space,  multiplicity mp$(M_\varphi)$ has the  meaning one expects intuitively, i.e., it equals the maximal number in $\N\cup\{\infty\}$ of preimages of $z\in\C$ under some $\varphi'=\varphi$ $\mu$-a.e. See (\ref{LMP}) for some details.

The  result we are going to prove is

\begin{The}\label{Main} Let $T$ be 
%a normal operator in a separable Hilbert space or 
a multiplication operator in $L^p(\mu)$, 
$p\in\,]0,\infty[$\,, on a finite separable measure space. Then  $\operatorname{mc}(T)=\operatorname{mc}^*(T)= \operatorname{mp}(T)$ holds. 
\end{The}

As already mentioned, due to the spectral theorem  the classical theorem of Bram \cite{Br}, by which any $*$-cyclic bounded normal operator is cyclic, is generalized by (\ref{Main}) to unbounded normal operators.  By definition  mp$(T)=1$ holds  if and only if  $T$ is isomorphic with $M_z$ in $L^p(\mu)$ for some finite Borel measure on $\C$. Hence (\ref{Main}) includes also the result that $M_z$ is cyclic.
Recall that a normal operator $T$ is said to be simple if its spectral measure is simple. So by  (\ref{Main})  $T$ is simple if and only if $T$ is cyclic. 
\\

If for a cyclic set $Z$ for $M_\varphi$ the subspace $\Pi(M_\varphi)Z$  is even a core of $M_\varphi$ 
then $Z$ is called \textbf{graph cyclic}. We have taken this expression from Szafraniec \cite{Sz}, which we consider appropriate in view of (\ref{IAeGZ}). In case that $1_\C$ is graph cyclic for $M_z$ in $L^2(\mu)$ then the Borel measure $\mu$ on $\C$ is called \textbf{ultradeterminate} by Fuglede \cite{Fu}.
One has

\begin{Pro}\label{IGCC} Let $p\in\,]0,\infty[$. Let $(\Omega,\mathcal{A}, \mu)$ be a finite separable measure space  and let $\varphi:\Omega\to \C$ be measurable.  If  $Z$ is a  cyclic  set for $M_\varphi$ in  $L^p(\mu)$, then
 $Z \text{e}^{-|\varphi|} $  is graph cyclic for $M_\varphi$.
\end{Pro}

In particular (\ref{IGCC}) shows that every cyclic normal operator is even graph cyclic. \\

For the proof of (\ref{Main}) we had first to establish that  $M_z$ in $L^p(\mu)$ is cyclic for every finite Borel measure $\mu$ on $\C$. More precisely we have 

\begin{The} \label{IMZC} Let $\mu$ be a finite Borel measure on $\C$. Then there is a positive Borel measurable function $\rho$ 
such that $\Pi(z)\rho$ is dense in $L^p(\mu')$ for all $p\in\,]0,\infty[$ if  
$\mu'$ is a  finite Borel measure on $\C$ equivalent to $\mu$. Moreover $h$ is  cyclic for $M_z$  in $L^p(\mu')$ if $h$ is  Borel measurable and  satisfies $0<|h|\le C\rho$ for some constant $C>0$. 
\end{The}

An immediate consequence of (\ref{IMZC})  due to Nagy  \cite{Na1} concerns polynomial approximation. It generalizes the result in Conway \cite[V.14.22]{Co1} for measures with compact support.  {\it Let $\mu$ be a Borel measure on $\C$ and let $f:\C\to\C$ be measurable. Then there is a sequence $(p_n)$ of  polynomials in $z$ with $p_n\to f$ $\mu$-a.e.} Indeed, without restriction $\mu$ is finite. Let $h:=\inf\{\rho,\frac{1}{1+|f|}\}$. Then $h$ is positive cyclic  and $fh$ is bounded.  Therefore there is a sequence $(p_n)$ satisfying $p_nh\to fh$ in $L^p(\mu)$, and the result follows  for some subsequence of $(p_n)$. \\

Cyclicity of $M_z$ in $L^2(\mu)$ has already been tackled by B\'ela Nagy in  \cite{Na1}  adapting in parts the original proof of  Joseph Bram \cite{Br} for bounded normal operators. See also Conway \cite[V.14.21]{Co1}  for a proof of Bram's theorem. The first  step (i) and important result  achieved  in \cite{Na1} is the decomposition of the complex plane into a null set and a countable union of increasing $\alpha$-sets. Secondly (ii) $\bar{z}$ is approximated by polynomials in $L^2(\sigma)$ for some finite Borel measure $\sigma$ on $\C$ equivalent to the original measure $\mu$. The third step (iii) in \cite{Na1} deals with the proof for the denseness  in $L^2(\sigma)$ of  the polynomials $\Pi(z)$. However the result obtained by Hilbert space methods is valid only for bounded  $M_z$. Indeed, \cite{Na1} starts  the third step with the (unfounded) assumption  that any function in the closure of $\Pi(z)$ is still square-integrable if multiplied by $z$. In other words, $\overline{\Pi(z)}\subset \mathcal{D}(M_z)$ is assumed. Proceeding on this assumption \cite{Na1}  shows $\overline{\Pi(z)}=L^2(\sigma)$ by a reducing subspace argument. Consequently   $\mathcal{D}(M_z)$ is  the whole of $L^2(\sigma)$ implying that $M_z$ is bounded, whence Nagy  \cite{Na1} does not achieve its goal. In addition, in accomplishing the reducing space argument, \cite{Na1} uses $*$-cyclicity of $M_z$ relying on a reference, which proves to be erroneous.
\\
\hspace*{6mm} Our first step (i)  in proving cyclicity of  $M_z$ in $L^p(\mu)$ for $p\in\,]0,\infty[$ and every finite Borel measure $\mu$ on $\C$   is the same as in \cite{Na1}. We present a short proof  (\ref{ICAPg}) of the decomposition valid for a large class of polish spaces including e.g. separable Banach spaces with real dimension $\ge 2$. In the second step  (ii) we show by induction that even $\Pi(z,\bar{z})$ is contained in the closure of $\Pi(z)$, see (\ref{IZZQ}). At this stage, in the third step (iii), we bring in $*$-cyclicity  of  $M_z$ in $L^2(\mu)$ by (\ref{IGSCV}) and thus  avoid a reducing subspace argument, which anyway is not available in the case $p\ne 2$.\\

In the sequel we denote by $\Pi(f_1,\dots,f_n)$ the set of complex polynomials in functions $f_1,\dots,f_n$ on some set with $f_i^{\,0}:=1$. 
Let $d\in\N$ and $|x|:=\sqrt{x_1^2+\dots +x_d^2}$ for $x=(x_1,\dots,x_d)\in\R^d$. Moreover, let $x_i$ also denote the $i$-th coordinate function on $\R^d$. Let $c>0$.

 \begin{Pro} \label{IGSCV} Let $p\in\,]0,\infty[$ and let $\mu$ be a finite Borel measure on $\R^d$.   Then  $\Pi(x_1,\dots,x_d)\operatorname{e}^{-c|x|^2}$ is
dense in $L^p(\mu)$. In particular, $\Pi(z,\bar{z})\operatorname{e}^{-c|z|^2}$ is
dense in $L^p(\mu)$ for any finite Borel measure $\mu$ on $\C$.
 \end{Pro}

 As an application of  mc$(T)= $ mc$^*(T)$ by (\ref{Main}) and of  (\ref{IGSCV}) we state that if $T$ is a  normal operator and if 
   $f\in \mathcal{D}(|T|^{2n})$ $\forall$ $n$, then $T$ is cyclic if and only if  e$^{-|T|^2}f$ is $*$-cyclic.\\

As already mentioned, (\ref{IGSCV}) is a corollary to the $*$-cyclicity of $M_z$ in $C_0(\C)$.   The question we pose  is about $(*)$-cyclicity of $M_z$ on $C_0(X)$ for $X\subset \C$.

\begin{The}  \label{CCOX}  
Let $X\subset \C$ be a  locally compact subspace. Then 
$M_z$ in $C_0(X)$ is $*$-cyclic by $e^{-c|z|^2}\eta$ with any positive $\eta\in C_0(X)$, and $M_z$ is cyclic if and only if every compact $K$ contained in  $X$ is an $\alpha$-set.
\end{The}

In view of (\ref{CCOX}) we remark that 
a locally compact subspace of $\C$ is $\sigma$-compact and locally closed. Hence,  if $X\subset \C$ is locally compact and every compact  $K\subset X$ has empty interior then  $X$ is nowhere dense. 
If $X\subset\C$ has empty interior and $\C\setminus K$ is connected for every compact $K\subset X$, then 
$\C\setminus X$ is dense and  has no bounded components, and vice versa.\\
 \hspace*{6mm} 
If  $X$ is compact  then $M_z$ in $C_0(X)$ is cyclic if and only if $1_\C$ is cyclic for $M_z$.  This is due to $||ph-fh||_{\infty,X}\ge C||p-f||_{\infty,X}$ with $C:=\inf_{z\in X}|h(z)|>0$ for $f\in C(X)$, $p\in\Pi(z)$, and $h$ a cyclic vector for $M_z$. Hence in case of compact $X$ one recovers Lavrentev's theorem on $\alpha$-sets from (\ref{CCOX}).\\
 \hspace*{6mm}  As an example, (\ref{CCOX}) implies that $M_z$ is cyclic by the function $\operatorname{e}^{-c|z|^2}$ in $C_0(X)$, where $X$ is the spiral $\{\operatorname{e}^{(1+i)t}:t\in \R\}$.\\
 \hspace*{6mm} 
 In this context we mention the  result by Lavrentev/Keldych \cite{St} that for a closed subset $X$ of $\C$ every continuous function on $X$ can be approximated uniformly by entire functions if and only if $\C\setminus X$ is dense, has no  bounded components, and is locally connected at infinity.

\section*{Proofs}

If necessary in order to avoid ambiguities we write $\overline{M}^\mu$ for the closure in  $L^p(\mu)$ of the subset $M$. Similarly $||f||_{p\mu}$ denotes the norm of $f\in L^p(\mu)$. We start   with two preliminary elementary results.

    \begin{Lem} \label{IACTD} Let $p\in\,]0,\infty[$. Let  
$(\Omega,\mathcal{A},\mu)$ be a measure space.  Let $M\cup \{h\}$ be a set of measurable functions on 
$\Omega$ with   $h\ne 0$ $\mu$-a.e. 
Then $Mh$ is dense in $L^p(\mu)$  if and only if $M$ is dense in $L^p(|h|^p\mu)$.
 \end{Lem} 
   
{\it Proof.} 
Set $\nu:=|h|^p\mu$.  --- Suppose $\overline{Mh}^{\,\mu} = L^p(\mu)$. Let  $g\in L^p(\nu)$, $\epsilon>0$. Then $gh\in L^p(\mu)$ and there is an $f\in M$ such that $\epsilon >||fh-gh||_{p\mu} =||f-g||_{p\nu}$. This proves $\overline{M}^{\,\nu}=L^p(\nu)$. --- Now suppose $\overline{M}^{\,\nu}=L^p(\nu)$ and let $f\in L^p(\mu)$, $\epsilon>0$. Then $f/h\in L^p(\nu)$ and there is $g\in M$ with $\epsilon>||f/h-g||_{p\nu} =||f-gh||_{p\mu}$. This proves  $\overline{Mh}^{\,\mu} = L^p(\mu)$.    \hfill{$\Box$}\\

    \begin{Lem} \label{IOSCVg} Let $p\in\,]0,\infty[$. Let  
$(\Omega,\mathcal{A},\mu)$ be a measure space.  Let $M\cup \{h\}$ be a set of measurable functions on 
$\Omega$ with $h$ bounded and   $h\ne 0$ $\mu$-a.e. If $M$ is dense in $L^p(\mu)$, then so is  $Mh$.
 \end{Lem} 

{\it Proof.} Let $C>0$ be a constant with $|h|<C$. Set $A_n:=\{|\frac{1}{h}|\le n\}$. For $\Delta \in \mathcal{A}$ with $\mu(\Delta)<\infty$ and
$\epsilon>0$ there exists $f\in M$  satisfying $||1_{\Delta \cap A_n} \frac{1}{h} - f||_p <\epsilon/C$. Then $||1_{\Delta \cap A_n}  - fh||_p <\epsilon$ holds, which implies $1_{\Delta \cap A_n}  \in \overline{Mh}$ for all $n\in\N$. Therefore $1_\Delta\in \overline{Mh}$ for  all $\Delta \in \mathcal{A}$ with $\mu(\Delta)<\infty$. The result follows.  \hfill{$\Box$}\\

As to the proof of  (\ref{ISCVFCO}) note that  $\Pi(x_1,\dots,x_d)\text{e}^{-c|x|^2}$ is not a subalgebra of $C_0(\R^d)$,  whence the Stone-Weierstra\ss\ theorem cannot be applied directly. In \cite{AF} a combination of the theorems of Hahn/Banach, Riesz, and Bochner is used to overcome this problem.

\begin{Pro} \label{ISCVFCO}
  $\Pi(x_1,\dots,x_d)\text{e}^{-c|x|^2}$ is dense in $C_0(\R^d)$. In particular,
 $\Pi(z,\bar{z})\text{e}^{-c|z|^2}$  is dense in $C_0(\C)$.
 \end{Pro}
 
 {\it Proof.}
 For convenience let $c=2$. The subalgebra
 $\Pi(x_1,\dots,x_d,\text{e}^{-|x|^2})\text{e}^{-2|x|^2}$ of $C_0(\R^d)$ satisfies the assumptions of the  Stone/Weierstra\ss\ Theorem.  Thus  
it is dense   in $C_0(\R^d)$. Therefore it remains to show $\Pi(x_1,\dots,x_d,\text{e}^{-|x|^2})\text{e}^{-2|x|^2} \subset \overline{\Pi(x_1,\dots,x_d)\text{exp}(-2|x|^2)}$, which follows from $\Pi(x_1,\dots,x_d)\text{e}^{-n|x|^2}  \text{e}^{-2|x|^2} \subset \overline{\Pi(x_1,\dots,x_d)\text{exp}(-2|x|^2)}$ for $n=0,1,2\dots$ by forming the linear span at the left hand side. Now this is shown by induction on $n$. For the step 
 $n \to n+1$ let $T_k$ denote the $k$-th Taylor polynomial of $\text{e}^z$ and let $p\in \Pi(x_1,\dots,x_d)$. Then
  $$
 || p\operatorname{e}^{-(n+1)|x|^2}\operatorname{e}^{-2|x|^2}
 - pT_k(-|x|^2)\operatorname{e}^{-n|x|^2}\operatorname{e}^{-2|x|^2}||_\infty
  \le C\, \,||\text{e}^{-|x|^2}\big(\operatorname{e}^{-|x|^2}-T_k(-|x|^2)\big)||_\infty$$  
with $C:=   || p\operatorname{e}^{-(n+1)|x|^2}||_\infty<\infty $.
Estimating the remainder function according to Lagrange one gets
  $\text{e}^{-t}| \text{e}^{-t} -T_k(-t)|= \text{e}^{-t}| R_{k+1}(-t,0)|=\text{e}^{-t}  \frac{ \text{e}^{\tau}}{(k+1)!}t^{k+1}  \le  \frac{ \text{e}^{-t}}{(k+1)!}t^{k+1}   $ with maximum at $t=k+1$, and  Sterling's formula yields   $\frac{ \text{e}^{-(k+1)}}{(k+1)!}(k+1)^{k+1}  \le \big(2\pi (k+1)\big)^{-1/2} \to 0$ for $k\to \infty$. \hfill{$\Box$}\\
  
Let $h_n$ denote the $n$-th Hermite function in one real variable. Then (\ref{ISCVFCO}) means that $h_{n_1}\times\dots\times h_{n_d}$, $n_1,\dots,n_d\in\N\cup\{0\}$ is total in $C_0(\R^d)$. --- In particular,
  for every continuous function $f$ on $\C$ vanishing at infinity there is a sequence $(p_n)$ of polynomials in $z$ and $\bar{z}$ such that $p_ne^{-|z|^2} \to f$ uniformly on $\C$.

\begin{PO} \textbf{(\ref{IGSCV}).}
Recall that $C_0(\R^d)$ is dense in $L^p(\mu)$  (see e.g. \cite{Cr}) and note that $||\cdot||_p\le \mu(\R^d)^{1/p}||\cdot||_\infty$.  Therefore the result follows from  (\ref{ISCVFCO}).\hfill{$\Box$}
\end{PO}

\begin{Lem} \label{ICAPg}  Let $X$ be a polish space where every pair of  distinct points are joint by infinitely many non-intersecting paths.  Let $\mu$ be a $\sigma$-finite Borel measure on $X$. Then there is an increasing sequence of  
 compact sets $F_n$  with empty interior and connected complement such that $\mu(X\setminus \bigcup_nF_n)=0$.

\end{Lem}
{\it Proof.}  Without restriction let  $\mu$ be  finite. 
All subspaces of $X$ are separable. Choose a countable dense set  $\{a_1,a_2,\dots\}$ in the complement of the set of mass points. Since the latter is countable, it does not contain an inner point by Baire's theorem, whence  $\{a_1,a_2,\dots\}$ is dense in $X$. Since $\mu$ is finite and since  there are  
infinitely many non-intersecting paths joining $a_n$ to $a_{n+1}$, for every $m\in\N$ there are connected measurable sets $A_n$ with $a_n,a_{n+1}\in A_n$ such that 
 $B_m:=\bigcup_n A_n$ is dense connected with $\mu(B_m) < \frac{1}{2m}$.
 ---  By Ulam's theorem (see, e.g.,  \cite[(2.67)] {Ca2}) $\mu$ is tight and in particular outer regular. Therefore
there is  an open $V_m$ with $B_m\subset  V_m$ and $\mu(V_m)<\frac{1}{m}$ and
there is an increasing sequence $(C_n)$  of compact sets with 
$\mu(X\setminus \bigcup_n C_n)=0$. \\
\hspace*{6mm} Now set $U_n:= \bigcup_{m\ge n}V_{2^{m}}$ and
$F_n:=C_n\setminus U_n$. Clearly,  $(F_n)$ is increasing, $\mu(\complement F_n)\to 0$, and $F_n$ is  compact. Its interior is empty since $B_{2^{n}}\subset U_n\subset \complement F_n$ is dense. $\complement F_n$ is connected. Indeed, let $U,V$ be open sets covering $\complement F_n$ with $U\ne \emptyset$ and  $U\cap V\cap \complement F_n=\emptyset$. 
%Note $A\subset \complement F_n$. 
Since $B_{2^{n}}$ is dense, $U\cap V=\emptyset$ follows, and since $B_{2^{n}}$ is connected, $V=\emptyset$ follows.
\hfill{$\Box$}\\

Note that separable Banach spaces with real dimension $\ge 2$ satisfy the assumptions on $X$ in (\ref{ICAPg}).   For $X=\C$ the proof can be further shortened taking in place of all $B_m$  a  single dense null set $B$ consisting of countably many straight lines through one common point, and of course, Ulam's theorem is not needed.  The first  (more cumbrous) proof for $X=\C$ and any finite Borel measure $\mu$ is given in Nagy \cite{Na1}. If  $\operatorname{supp}(\mu)$ is compact   there is the original proof by Bram \cite{Br}, a similar one in Conway \cite{Co1}, and a simpler one in Shields \cite{Sh}.  \\

\begin{Lem} \label{IZZQ}  Let  $\mu$ be a finite   Borel measure on   $\C$. Then there is a positive Borel measurable function $\rho$ %with $0<h\le \operatorname{e}^{-c|z|^2}$ 
such that 
$\Pi(z)\subset L^p(\nu)$  and $\Pi(z,\bar{z})\subset \overline{\Pi(z)}^{\,\nu}$ for  all $p\in\,]0,\infty[$  and all  $\nu:=|h|^p\mu'$ with  
 $\mu'$ a finite Borel measure on   $\C$ equivalent  to $\mu$ and $h$ Borel measurable satisfying $0<|h|\le \rho$.

\end{Lem}

{\it Proof.} 
Set $k:\C\to\C$, $k(z):=\bar{z}$.
 Let $(F_n)$ be an increasing sequence of $\alpha$-sets  of $\C$ satisfying $\mu(N)=0$ for $N:=\C\setminus \bigcup_nF_n$, see (\ref{ICAPg}). 
 For every $n\in\N$ there is  $q_n\in\Pi(z)$ satisfying $||1_{F_n}(k-q_n)||_\infty<\text{e}^{-\delta_n}$ with $\delta_n:=n || 1_{F_n} k||_\infty$. Set $M_n:=\max\{1,||q_1\text{e}^{-|z|}||_\infty,\dots,||q_n\text{e}^{-|z|}||_\infty\}$ and let $\rho$ be the positive function on $\C$ given by $\rho|N:=1$ and $\rho|(F_n\setminus F_{n-1}):=\text{e}^{-2|z|}/M_n$ for $n\in\N$ with $F_0:=\emptyset$. \\%Let $\nu:=\rho^p\mu$.\\
\hspace*{6mm} Since $q\rho$  for $q\in\Pi(z)$ is bounded, $\Pi(z)\subset L^p(\nu)$ holds.  For the proof of $\Pi(z,\bar{z})\subset \overline{\Pi(z)}^{\,\nu}$ obviously  it suffices to show $\Pi(z)\bar{z}^m\subset \overline{\Pi(z)}^{\,\nu}$  for $m=0,1,2,\dots$ 
This occurs by induction on $m\in\N\cup\{0\}$. Let $j\in\N\cup\{0\}$ and write $\bar{z}^{m+1}=\bar{z}^mk$. Then 
$|| z^j \bar{z}^m  k-  z^j \bar{z}^m q_n||_{p\nu}\le \nu(\C)^{1/p}\delta_n^{(j+m)}\,  \text{e}^{-\delta_n}+||1_{\complement F_n}(z^j\bar{z}^mk-z^j\bar{z}^mq_n)||_{p\nu}$. The first summand vanishes for $n\to\infty$,
 the latter is less or equal to $||1_{\complement F_n}z^j\bar{z}^mk||_{p\nu} + ||1_{\complement F_n}z^j\bar{z}^mq_n\rho||_{p\mu'}$, up to the constant factor $\sqrt[p]{2}/2$ %$2^{\frac{1}{p}-1}$ 
in the case $p<1$. Now $|1_{\complement F_n} z^j\bar{z}^mk|\le|z^{(j+m+1)}|$ and $|1_{\complement F_n}z^j\bar{z}^mq_n\rho|\le 
|z^{(j+m)}|\text{e}^{-|z|}$, whence both summands vanish for $n\to \infty$ by dominated convergence. Since 
$z^jq_n\bar{z}^m\in \overline{\Pi(z)}^{\,\nu}$ by assumption,
we infer $z^j\bar{z}^m k\in  \overline{\Pi(z,\bar{z})}^{\,\nu}$ for every $j$, thus concluding the proof.\hfill{$\Box$}\\

\begin{PO} \textbf{(\ref{IMZC}).} By (\ref{IZZQ}) there exists a Borel measurable function $\rho$  with $0<\rho\le \text{e}^{-|z|^2}$ such that $\Pi(z)\subset L^p(\nu)$  and $\Pi(z,\bar{z})\subset \overline{\Pi(z)}^{\,\nu}$ for  all $p\in\,]0,\infty[$  if   $\nu:=\rho^p\mu'$.  
By (\ref{IGSCV}) and (\ref{IOSCVg})  we have $\overline{\Pi(z,\bar{z})\rho}^{\,\mu'}=L^p(\mu')$, whence $\overline{\Pi(z,\bar{z})}^{\,\nu}=L^p(\nu)$ by (\ref{IACTD}). Therefore $\overline{\Pi(z)}^{\,\nu}=L^p(\nu)$ and hence $\overline{\Pi(z)\rho}^{\,\mu'}=L^p(\mu')$ by (\ref{IACTD}). The last assertion  follows from  (\ref{IOSCVg}) for $M=\Pi(z)$. 
\end{PO}

The decomposition $(\pi,\nu)$ of $(\varphi,\mu)$ can be derived from Rohlin's disintegration theorem (see Rohlin \cite{Ro}), and  can be found in  Seid  \cite[remark 3.4]{Se}. There $\varphi$ is a bounded Borel function on the finite measure space $([0,1], \mathcal{B}, \mu)$ with $\mathcal{B}$ the Borel sets. By the following two lemmata, which we state without proof,  we generalize this result in that $\varphi$ is a measurable  not necessarily bounded function on a finite separable measure space $(\Omega,\mathcal{A},\mu)$  and the multiplication operator $M_\varphi$ is isomorphic with $M_\pi$ by means of a measure algebra isomorphism.

\begin{Lem} \label{FL} Let $(\Omega,\mathcal{A},\mu)$ be a finite  separable measure space. Then there is an $(\mathcal{A},\mathcal{B})$-measurable function $a: \Omega\to [0,1]$ such that 
$[B] \mapsto [a^{-1}(B)]$, $B\in \mathcal{B}$ is a measure algebra isomorphism from  $(\mathcal{B},a(\mu))$ onto $(\mathcal{A},\mu)$ and that 
for every measurable $\varphi:\Omega\to \C$ there is an $a(\mu)$-almost unique  measurable $\psi: [0,1] \to \C$ with $\varphi= \psi \circ a$ $\mu$-a.e.  If $(\Omega,\mathcal{A})$ is a standard measurable space then a measurable space isomorphism $a$ exists.
\end{Lem}

\begin{Lem} \label{UBP} Let $(\Omega,\mathcal{A},\mu)$ be a finite separable measure space and let $\varphi:\Omega \to \C$ be measurable. Let $k:\C\to \mathbb{D}$, $k(z):=\frac{z}{1+|z|}$, and let $\kappa$ denote its inverse. Let $\gamma: [0,1] \times \mathbb{D} \to [0,1] \times \C$, $\gamma(t,z):=(t,\kappa(z))$. Then $k\circ \varphi$ maps to $\mathbb{D} $ and, if $(\pi,\nu)$ is a Rohlin decomposition of $(k\circ \varphi,\mu)$, then $(\pi,\gamma(\nu))$ is a Rohlin decomposition of $(\varphi,\mu)$.
\end{Lem}

\begin{PO} \textbf{(\ref{Main}).} By definition of multicyclicity, multi-$*$-cyclicity, and  multiplicity it suffices to show   $\operatorname{mc}(M_\pi)=\operatorname{mc*}(M_\pi)=\operatorname{mp}(M_\pi)$ for a Rohlin decomposition $(\pi,\nu)$. Plainly
$\operatorname{mc*}(T)\le\operatorname{mc}(T)$.\\
\hspace*{6mm} Consider first the case $\mu_c\ne 0$. Then $\operatorname{mp}(M_\pi)=\infty$ holds by definition. In oder to show  $\operatorname{mc*}(M_\pi)=\infty$, it suffices to treat the case $\nu=\lambda\otimes \mu_c$,
since $L^p(\nu)$ is the direct sum of the subspaces $L^p(\lambda\otimes\mu_c)$ and $L^p(\delta_{1/n}\otimes \mu_n)$, which are invariant under $M_\varphi$ and $M_{\overline{\varphi}}$. Let us assume that there is a finite $*$-cyclic set  $\{f_1,\dots,f_d\}$ for $M_\pi$ with $d\in\N$. Set $\chi_n:=1_{]1/(n+1),1/n]}$.  For $\chi_n1_\C$ there are sequences $(p_{n\delta k})_k$ in $\Pi(z,\bar{z})$ such that $q_{n\delta k}(z):=p_{n\delta k}(z,\bar{z})$ satisfy $\sum_{\delta=1}^d q_{n\delta k}f_\delta\to \chi_n1_\C$ f\"ur $k\to\infty$ in $L^p(\nu)$. Set $f_{\delta z}:=f_\delta(\cdot,z)$. By Tonelli's theorem there is a subsequence $(k_l)_l$, without restriction $(k)_k$ itself, with 
 $\sum_{\delta=1}^d q_{n\delta k}(z) f_{\delta z}    \to \chi_n$ for $k\to\infty$ in $L^p(\lambda)$ for $\mu_c$-a.a. $z\in\C$. We consider this convergence for $n=1,\dots,d$.  
Since  $\chi_1,\dots,\chi_d$ are  linear independent, it follows that $f_{1 z},\dots,f_{d\, z}$ are so  for  $\mu_c$-a.a. $z\in\C$. 
Consequently $q_{n\delta k}(z)$ for $k\to \infty$ converge to the coordinates $\alpha_{\delta n}(z)$ of  $\chi_n$ with respect to $(f_{\delta z})_\delta$.   Hence one gets  $f_{\delta z}=\sum_{n=1}^d\beta_{n \delta}(z)   \chi_n$, $\delta=1,\dots,d$ with coordinates  $\beta_{n \delta}(z)$.  --- Repeating these considerations for $\chi_{d+1},\dots,\chi_{2d}$ in place  of $\chi_1,\dots,\chi_d$ we obtain $f_{\delta z}=\sum_{n=1}^d\beta'_{n \delta}(z)   \chi_{d+n}$, $\delta=1,\dots,d$. This implies  for 
$\mu_c$-a.a. $z\in\C$ that $ f_\delta(t,z)=0$ for  $\lambda$-a.a. $t\in[0,1]$. Hence   $\int |f_\delta|^p\text{d}\, \lambda\otimes\mu_c =0$ by Tonelli's theorem. This means   $f_\delta =0$ for $\delta=1,\dots,d$, which is not possible.\\
\hspace*{6mm}  Now let $\mu_c=0$ and set $\chi_n:=1_{\{1/n\}}$, $n\in\N$. We consider first the case 
$\operatorname{mp}(M_\pi)=\infty$. Then $\mu_n\ne 0$ $\forall$ $n\in\N$. Assuming the existence of a $*$-cyclic set of $d$ elements, as in the previous case $f_{\delta z}=\sum_{n=1}^d\beta_{n \delta}(z)   \chi_n=\sum_{n=1}^d\beta'_{n \delta}(z)   \chi_{d+n}$, $\delta=1,\dots,n$ follows. This implies for $\mu_c$-a.a. $z\in\C$ that 
$f_\delta(t,z)=0$ for all  $1/t\in\N$. This means  $f_\delta=0$ for all $\delta=1,\dots,d$, which is not possible.  --- 
We turn to the last case  $N:=\operatorname{mp}(M_\pi)\in\N$. Then $\mu_n\ne 0$ f\"ur $n=1,\dots,N$ und $\mu_n=0$ else.
Since $M_\pi$ is cyclic in $L^p(\delta_{1/n}\otimes \mu_n)$ according to (\ref{IMZC}), $\operatorname{mc}(M_\pi)\le N$ follows.
--- Let us assume now that there is a $*$-cyclic set  $\{f_1,\dots,f_d\}$ for $M_\pi$ with $d< N$. By   considerations as in the case 
$\mu_c\ne 0$ we get $f_{\delta z}=\sum_{n=1}^d\beta_{n \delta}(z)   \chi_n=\sum_{n=2}^{d}\beta'_{n \delta}(z)   \chi_{n} +\beta'_{m\delta}\chi_m$ for some $m\not\in\{1,\dots,d\}$. Since all $\chi$'s  are linear independent this implies 
$\beta_{1\delta}(z)=0$ f\"ur $\mu_c$-a.a. $z\in\C$. Analogously $\beta_{n\delta}(z)=0$ for $\mu_c$-a.a. $z\in\C$ for every $n\in\{1,,\dots,d\}$. Therefore $f_{\delta z}=0$ for $\mu_c$-a.a. $z\in\C$. This means $f_\delta=0$ for every  $\delta$, which is not possible.
\hfill{$\Box$}
\end{PO}

The   next lemma is not new but  it puts together the equivalences  for convenience.

\begin{Lem} \label{IAeGZ} Let $p\in\,]0,\infty[$. Let $(\Omega,\mathcal{A}, \mu)$ be a finite measure space and let  $Z\cup \{\varphi\}$ be a set of  measurable functions. Then \emph{(a) -- (d)} are equivalent.
\begin{itemize}
\item[\emph{(a)}] $Z$ is graph cyclic for $M_\varphi$ in $L^p(\mu)$
\item[\emph{(b)}]  $\Pi(M_\varphi)Z$ is a core for $M_\varphi$ in $L^p(\mu)$
\item[\emph{(c)}] $\{(f,M_\varphi f): f \in \Pi(M_\varphi)Z \}$ is dense in $\{(f,M_\varphi f): f \in \mathcal{D}(M_\varphi)  \}$ 
\item[\emph{(d)}] $\Pi(M_\varphi)Z \sqrt[p]{1+|\varphi|^p}$ is dense in $L^p(\mu)$
\item[\emph{(e)}]  $\Pi(M_\varphi)Z$ is dense in $L^p\big((1+|\varphi|^p)\,\mu\big)$
\end{itemize}
\end{Lem} 

{\it Proof.} The equivalences of (a) and (b)  and (c) hold by definition,  the equivalence of (d) and (e) holds by (\ref{IACTD}).  --- As to (c) $\Rightarrow$ (d) let $g\in L^p(\mu)$. Then $g':=g\big/\sqrt[p]{1+|\varphi|^p} \in  \mathcal{D}(M_\varphi)$ and hence for $\epsilon>0$ there is $f\in \Pi(M_\varphi)Z$ satisfying $|| (f,M_\varphi f)- (g',M_\varphi g')||_p <\epsilon$, which means $\epsilon^p>\int \big(|f-g'|^p+ |\varphi f- \varphi g'|^p\big)\operatorname{d}\mu=  \int \big|  \sqrt[p]{1+|\varphi|^p} \, f-g\big|^p \operatorname{d}\mu$ proving (d). --- Finally assume (d) and let $ f \in \mathcal{D}(M_\varphi) $. Then $f':= \sqrt[p]{1+|\varphi|^p}\, f \in L^p(\mu)$
  and  for $\epsilon>0$ there is $g\in \Pi(M_\varphi)Z$  satisfying $||f' - \sqrt[p]{1+|\varphi|^p}\,g||_p<\epsilon$, which means $|| (f,M_\varphi f)- (g,M_\varphi g)||_p <\epsilon$, thus proving (c). \hfill{$\Box$}\\

\begin{PO} \textbf{(\ref{IGCC}).}  Since $\Pi(M_\varphi)Z$ is dense in $L^p(\mu)$, by (\ref{IOSCVg}) also $\Pi(M_\varphi)Z\text{e}^{-|\varphi|}\sqrt[p]{1+|\varphi|^p}$  is dense in $L^p(\mu)$. The result follows from (\ref{IAeGZ}).      \hfill{$\Box$}

\end{PO}

\begin{PO} \textbf{(\ref{CCOX}).} Let $c=1$ for convenience. ---  (i) Since $X$ is locally compact there are compact sets $F_n\subset X$ with $F_n$ contained in the interior of $F_{n+1} $, and functions $\eta_n\in C_c(X)$ satisfying $\eta_n|F_n=1$, $\eta_n|(X\setminus F_{n+1})=0$, and $0\le \eta_n\le 1_X$. Let $\alpha_n>0$ with $\sum_n\alpha_n<\infty$. Then $\eta:=\sum_n\alpha_n\eta_n\in C_0(X)$ and $\eta>0$.\\
\hspace*{6mm} (ii) Let  $f \in C_c(X)$. Extend $f$ continuously, first onto the closure $\overline{X}$ by $0$, and subsequently onto $\C$ by the Tietze-Urysohn extension theorem. Finally, multiplying the resulting function  by a $j\in C_c(\C)$ with $j|\operatorname{supp}(f)=1$ one achieves an extension of $f$ with compact support.
 --- Now let $g\in C_c(\C)$  extend  $f/\eta \in C_c(X)$ with a positive $\eta \in C_0(X)$, see (i).
 Let $\epsilon >0$. By (\ref{ISCVFCO})
there is $p\in\Pi(z,\bar{z})$ with $ ||g-p\operatorname{e}^{-|z|^2}||_\infty <\epsilon/||\eta||_{\infty,X}$. This implies 
$||f- p\operatorname{e}^{-|z|^2}\eta||_{\infty,X} <\epsilon$. Thus $e^{-|z|^2}\eta$ is a $*$-cyclic vector.\\
\hspace*{6mm} (iii) Let $h\in C_0(X)$ be a cyclic vector. Since $C_0(X)$ vanishes nowhere, so does $h$. Let $K\subset X$  be  compact. Let $\varphi\in C(K)$. By the Tietze-Urysohn extension theorem exists a bounded continuous $\phi$ on $X$ with $\phi|K=\varphi$. Then $\phi h\in C_0(X)$. Set $c:= \sup_{z\in K} |\frac{1}{h(z)}|$ and let $\epsilon>0$.  Then there is $p\in\Pi(z)$ satisfying $||\phi h-ph||_\infty<\epsilon/c$. This implies $||\varphi-p||_{\infty, K}<\epsilon$. Thus $K$ is an $\alpha$-set by definition.\\ 
\hspace*{6mm} (iv) Now let every compact $K\subset X$ be an $\alpha$-set. Set $k(z):=\bar{z}$. There are $q_n\in \Pi(z)$ satisfying $||1_{F_n}(k-q_n)||_\infty<\frac{1}{n}$. Then set 
$M_n:=\max\{1,||1_{F_{n+1}}k||_\infty,||1_{F_{n+1}}q_1||_\infty,\dots,$
$||1_{F_{n+1}}q_n||_\infty\}$ and $\alpha_n:=2^{-n}/M_n$ in (i).
For $j\ge n$ one has $||\alpha_j\eta_j(k-q_n)||_{\infty,X}\le 2\cdot 2^{-j}$, whence $||1_{\complement F_n}\eta (k-q_n)||_{\infty,X} \to 0$. It follows $||\eta (k-q_n)||_{\infty,X} \to 0$. ---  Now we show that $h:=\operatorname{e}^{-|z|^2}\eta$ is a cyclic vector. By (ii), $A:=\Pi(z,\bar{z})h$ is dense in $C_0(X)$.
We conclude the proof showing $A\subset \overline{\Pi(z)h}$ by the method used in (\ref{ISCVFCO}). Let $q\in\Pi(z)$. 
 Induction occurs on $m=0,1,2\dots$ Then
$||q\,\bar{z}^{m+1}h\,-\,q\,q_n\bar{z}^mh||_{\infty,X}\le C\;||\eta(k-q_n)||_{\infty,X}$
with $C:=||q\,\bar{z}^m\operatorname{e}^{-|z|^2}||_{\infty}$ vanishes for $n\to \infty$. \hfill{$\Box$}\\
\end{PO}

\section*{Further results}

Let $p\in\,]0,\infty[$. We know by (\ref{IGSCV}), (\ref{IMZC}) that, for every finite Borel measure $\mu$, $M_z$ in $L^p(\mu)$ is $*$-cyclic by the continuous vector  e$^{-c|z|^2}$ and that $M_z$ is cyclic.
The question is about continuity of a cyclic vector.

\begin{Exa}  \label{NCCV}  Let $\mu:=1_\mathbb{D}\lambda$ with $\lambda$ the Lebesgue measure on $\C$ and $\mathbb{D}$ the open unit disc. Let $h$ be a cyclic vector for $M_z$ in $L^2(\mu)$. Then $\{h= 0\}$ is a $\mu$-null set containing all continuity points of $h$. 
\end{Exa}

{\it Proof.}  $\{h=0\}$ is a $\mu$-null set, since $L^2(\mu)= \overline{\Pi(z)h}\subset \{f\in L^2(\mu): f=1_{\{h\ne0\}}f\}$. ---  Let $h$ be continuous at $x\in \mathbb{D}$. Assume $h(x)\ne 0$. Then there are an open disc $D$ with center  $x$ and $\delta>0$ such that $\delta1_D\le |h|$. Hence, by (\ref{IOSCVg}),  $\Pi(z) 1_D$ is dense in $L^2(1_D\lambda)$. This contradicts
e.g. 3.22. (c) in  \cite{Ca}.
\hfill{$\Box$}\\

In particular, there is no  cyclic vector for $ M_z$ in $L^2(1_\mathbb{D}\lambda)$ that is continuous on $\mathbb{D}$, thus answering  a question about continuity of cyclic vectors posed by Shields \cite{Sh}. If, however, $\bar{z} \in \overline{\Pi(z)}$ holds then we have

\begin{Pro}  \label{ICCV}  Let $p\in\,]0,\infty[$.  Let $\mu$ be a finite Borel measure on $\C$ such that $\Pi(z)\subset L^p(\mu)$ and $\bar{z} \in\overline{\Pi(z)}$, then   $\operatorname{e}^{-c|z|^2}$ for $c>0$ is a cyclic vector for $M_z$ in $L^p(\mu)$.

\end{Pro}
 {\it Proof.} Apply (\ref{tool}) to $A:=\Pi(z)$, $b:=\bar{z}$, and $c:=\operatorname{e}^{-c|z|^2}$. The result follows from (\ref{IGSCV}).
 \hfill{$\Box$}\\

The following is a useful tool in establishing as in (\ref{ICCV}) that the closure of a coset of a given algebra contains the coset of some larger algebra.

 \begin{Lem}\label{tool} Let $p\in\,]0,\infty[$ and let $\mu$ be a finite Borel measure on $\C$. Let $A \subset L^p(\mu)$ be an algebra, let $b\in \overline{A}$, and let $c\in L^p(\mu)$. Suppose  $Ab^nc\subset L^\infty(\mu)$ for $n=0,1,2\dots$. Then $\Pi(A,b)c \subset\overline{Ac}$.
 \end{Lem}
 
{\it Proof.}  It suffices to show $Ab^nc\subset \overline{Ac}$ by induction on $n$. Let $(a_k)$ be a sequence in $A$ converging to $b$. As to the step $n\to n+1$ note $||ab^{n+1}c-aa_kb^nc||_p\le ||ab^nc||_\infty ||b-a_k||_p \to 0$ for $k\to \infty$. Since by assumption $aa_kb^nc \in  \overline{Ac}$ the result follows.\hfill{$\Box$}\\

For $z\in\C$ and  $\Delta\subset \Omega$   let  n$_\Delta(z):=|\{\omega\in \Delta:\varphi(\omega)=z\}|$ denote the number in $\N\cup\{\infty\}$
of preimages   in $\Delta$ of $z$ under $\varphi$. For  a   Rohlin decomposition  $(\pi,\nu)$ of $(\varphi,\mu)$, $\varphi(\mu)=\mu_c+\sum_n\mu_n$ holds.  Set  $P_n:=\{\frac{\operatorname{d}\mu_n}{\operatorname{d}\varphi(\mu)}>0\}$ for $n\in\{c\}\cup \N$ and define the \textbf{local multiplicity} by
$$\operatorname{m}_\varphi(z):=\infty1_{P_c}(z)+\sup\,\{0\}  \cup\{n\in\N:\,z\in P_n\}.$$ 
We will keep in mind that $P_n$ is unique up to a $\varphi(\mu)$-null set.

 \begin{Lem}\label{LMP} Let $(\Omega,\mathcal{A},\mu)$ be a finite measure space with $(\Omega,\mathcal{A})$ a standard measurable space. Let $\varphi:\Omega\to \C$ be measurable. Let $(\pi,\nu)$ be a Rohlin decomposition of $(\varphi,\mu)$ by a point isomorphism  $\vartheta$ onto the complement of a $\nu$-null set of $[0,1]\times\C$. 
 Then there is a $\mu$-null set $N$ such that $\operatorname{n}_{\Omega\setminus N}$ is measurable and such that 
  $\operatorname{m}_\varphi=\operatorname{n}_{\Omega\setminus N}\le \operatorname{n}_{\Omega\setminus N'} $ $\varphi(\mu)$-a.e. for every $\mu$-null set $N'$. Furthermore  $\operatorname{mp}(M_\varphi)=\sup \operatorname{m}_\varphi$ holds  for $M_\varphi$ in $L^p(\mu)$, $p\in\;]0,\infty[$.
 \end{Lem}
 
 {\it Proof.} Because of $\mu_{n+1}\ll \mu_n$ $\forall$ $n$, $(P_n)_n$ is $\varphi(\mu)$-almost decreasing. Therefore m$_\varphi=\infty 1_{P_c}+\sum_n1_{P_n}$ $\varphi(\mu)$-a.e. holds.
 $S:=\big([0,1]\times P_c\big)\cup\bigcup_n\big(\{\frac{1}{n}\}\times P_n\big)$ is the complement of a $\nu$-null set since $\mu_n(\complement P_n)=0$ for $n=c,1,2,\dots$, and  $|S_z|=\infty 1_{P_c}(z)+\sum_n1_{P_n}(z)$ $\forall$ $z\in\C$ for $S_z:=\{t\in [0,1]:(t,z)\in S\}$ holds. --- Now let $R$ be the complement of any $\nu$-null set. Then $B_c:=\{z\in\C:\lambda(R_z)=1\}$ and $B_n:=\{z\in\C:\frac{1}{n}\in R_z\}$ satisfy $\mu_n(\C\setminus B_n)=0$, whence $\varphi(\mu)(P_n\setminus B_n)=0$ for $n=c,1,2,\dots$ This implies $|R_z|\ge |S_z|$ 
 $\varphi(\mu)$-a.e. Moreover,  we may choose without restriction  $P_n\subset B_n$, $n=c,1,2,\dots$  for $R:=\vartheta(\Omega)$.   Then $|S_z|=|S_z \cap \vartheta(\Omega)|$. Since generally $\operatorname{n}_\Delta(z)=\vartheta(\Delta)_z$ because of $\varphi=\pi \circ \vartheta$, we obtain  
$\operatorname{n}_{\Omega\setminus N}(z)=|S_z|$ for $N:=\complement \,\vartheta^{-1}(S)$, whence $\operatorname{n}_{\Omega\setminus N}$ is measurable, and $\operatorname{n}_{\Omega\setminus N}(z)=|S_z|\le |\vartheta(\Omega\setminus N')_z|=\operatorname{n}_{\Omega\setminus N'}(z)$ $\varphi(\mu)$-a.e. --- The last assertion is obvious.\hfill{$\Box$}\\

Obviously, in (\ref{LMP}), $\operatorname{n}_{\Omega\setminus (N\cup N')}= \operatorname{n}_{\Omega\setminus N}$, whence
 $\operatorname{m}_\varphi=\operatorname{n}_{\Omega\setminus M}$ holds $\varphi(\mu)$-a.e., if the $\nu$-null set $M$ is large enough. Finally we mention that in the Hilbert space case m$_\varphi$ is a complete invariant. This means that normal operators $T\simeq M_\varphi$ in $L^2(\mu)$ and $T'\simeq M_{\varphi'}$ in $L^2(\mu')$  with $\mu$ and $\mu'$   Borel measures  on $\C$ are isomorphic if and only if  $\varphi(\mu)\sim\varphi'(\mu')$ and  $m_\varphi=m_{\varphi'}$ a.e. In other words $m_\varphi$ is the usual local multiplicity derived from the spectral theorem. Results relating local multiplicity to the number of preimages can be found in \cite{Az,Ho,Ka,Kr,Nad}.

\end{document}